\documentclass[twoside,reqno,11pt]{article} %
   \usepackage{hyperref}
   \def\theequation{\arabic{section}.\arabic{equation}}

\usepackage{amsmath,amssymb,upgreek,enumerate}
\newtheorem{theorem}{Theorem}[section]
\newtheorem{corollary}[theorem]{Corollary}
\newcommand{\e}{\mathrm{e}}
\newcommand{\dd}{\mathrm{d}}
\newcommand{\D}{\mathrm{D}}
\newcommand{\I}{\mathrm{I}}

\newcommand{\ii}{\mathrm{i}}
\newcommand{\upi}{\mathrm{\pi}}

 \title{A comment on a controversial issue: \\ [3pt]
a Generalized Fractional Derivative cannot \\ [3pt] have a regular kernel}

\date{February 16, 2020}

\author{Andrzej Hanyga}

\begin{document}

\maketitle 

\begin{abstract}

The problem whether a given pair of functions can be used as the kernels of a generalized fractional derivative
and the associated generalized fractional integral
is reduced to the problem of existence of a solution to the Sonine equation. It is shown for some selected classes of
functions that a necessary condition
for a function to be the kernel of a fractional derivative is an integrable singularity at 0.
It is shown that locally integrable completely monotone functions satisfy the Sonine equation if and only if they are singular at 0.

 \medskip


 \smallskip


\end{abstract}

\textbf{Notations.}\\ 
%\vskip 0.1cm

\begin{tabular}{ll}
$[a,b[$ = $\{x \mid a \leq x < b\}$, \\
$f \ast g$ = $\int_0^t f(s)\, g(t-s) \, \dd s$,\\
$\theta^{\alpha}(t)  = t^\alpha/\Gamma(1 + \alpha)$.
\end{tabular}

%\vspace*{-25pt}

\section{Introduction} %%%%%%%%%%% Section 1 %%%%%%%%%%%%%%%%%%

\medskip

A few generalized fractional derivatives (GFD) with regular kernels have been recently recently suggested beginning with the Caputo-Fabrizio and the Atangana-Baleanu derivatives \cite{CF,AB}.
They were followed by an avalanche of applications in various disciplines although critical assessments of their applicability are lacking.
These papers can be found by scholar.google search using the keywords ``fractional derivative" and ``regular" or in the bibliography of \cite{QRB}.

Several paradoxes involving the GFD with regular kernels have been pointed out in \cite{St,OT1,T,G}. In \cite{OT2} it has been shown that
the models involving the GFD with regular kernels poorly reflect the real world data. All these paradoxes stem
from the fact the
GFD with regular kernels are equivalent to Volterra integral operators, while the associated generalized fractional integrals are not
integral operators. The generalized fractional integral of the Caputo-Fabrizio \cite{CF} or Atangana-Baleanu type \cite{AB} is a weighted sum of an identity operator and an integral, with the weights determined by
a parameter $\alpha \in [0,1[$, misleadingly called order. For $\alpha \rightarrow 0$ the generalized fractional integral of the Caputo-Fabrizio or Atangana-Baleanu type tends to the identity operator,
while for $\alpha \rightarrow 1$ it tends to a pure integral operator.

In several papers the GFD with regular kernels were substituted for ordinary time derivatives in several evolution equations.
Stynes \cite{St} showed that this procedure leads to ill-posed initial-value problems. Here is his argument in a simplified form.
 A Caputo-type GFD with a kernel bounded in a right neighborhood $[0,\varepsilon]$, $\varepsilon > 0$, of 0 has the form $\D_g \, \varphi(t) = g\ast \varphi^\prime$ and $\vert g \vert \leq K < \infty$. If $\varphi^\prime\vert_{[0,\varepsilon]} \in \mathcal{L}^1([0,\varepsilon])$, then
$$\vert \D_g \varphi(t) \vert \leq K \left\vert \int_0^t \varphi^\prime(s)\, \dd s\right\vert \rightarrow 0$$ for $t \rightarrow 0$. Hence an initial-value problem $$\D_g \, u(t,x) = L_x u(t,x) + k(t,x),\;\; u(0,x) = h(x),$$ where $L_x$ is an operator involving spatial derivatives, does not have solutions unless
$L_x \, h(x) + k(0,x) = 0$.

In the case of the Caputo-Fabrizio derivative \cite{CF} the derivative $g^\prime$ is locally integrable and $\D_g \, \varphi(t)
\equiv  g(0) \, \varphi(t) + g^\prime \ast \varphi$ with $g(0) > 0$. The Caputo-Fabrizio fractional equation $\D_g \, \varphi = F(t)$
is thus equivalent to a Volterra integral equation of the second kind and therefore cannot be interpreted as an evolution equation.

Specific fractional-differential equations were considered by Giusti \cite{G}, who showed that they are equivalent to differential equations of integer order or fractional order equations with Caputo fractional derivatives.

The objections raised in the critical papers cited above concern the shortcomings of the Caputo-Fabrizio and Atangana-Baleanu derivatives.
In contrast to the above-mentioned papers we investigate here the conditions to be met by the kernel of a generalized fractional derivative.
Our definition of a GFD assumes that the kernels of the derivative and the associated integral are
both locally integrable functions. It follows that the two kernels should satisfy the Sonine equation, defined below. It will be
concluded from this fact that under an additional assumption the kernels must be singular at 0.

The proof of the singularity of the solutions of the Sonine equation in the class of locally integrable functions requires a technical assumption
which is difficult to verify without knowing both functions (Section~\ref{sec:necessary}). The situation is much easier if one
of the functions is a priori known to be locally integrable and completely monotone (LICM). Every
singular LICM function is a Sonine function (Section~{\ref{sec:CM}). In this case  both associated Sonine functions are singular LICM.

The kernels of fractional derivatives and integrals as well as the kernels of the Caputo-Fabrizio, Atangana-Baleanu derivatives and other fractional derivatives with regular kernels put forward so far are LICM. Consequently as for now the results of Section~{\ref{sec:CM}} cover all the
candidates for generalized fractional derivatives discussed in the literature.

The singularity of a sLICM function can be logarithmic. A simple example of an sLICM function with a logarithmic singularity is
presented in Section~\ref{log}. 

It may be still a complicated task to find an explicit formula for the associated Sonine function for a given singular LICM function.
In the last section we show how the theory of completely monotone functions may be used to find the associated Sonine function for a given singular LICM function (Section~{\ref{sec:exist}}).

\vspace*{-3pt} 
\section{What is a generalized fractional derivative?} 

\setcounter{section}{2} \setcounter{equation}{0}

In this paper we shall consider a generalization of Riemann-Liouville derivatives of order $< 1$. This generalization includes
some operators of practical interest, viz. distributed-order fractional derivatives obtained by averaging fractional derivatives over orders $ \alpha \in [0,1]$
\cite{C,H3} and shifted fractional derivatives \cite{victoria}. The object of our study, which is the relation between the kernels of the generalized fractional derivatives and
generalized fractional integral and its implications, is identical in the case of Caputo-type generalized fractional derivatives
studied in \cite{K}.

A generalized fractional derivative $\D_g$ (in the Riemann-Liouville sense) and an associated generalized fractional integral
$\I_f$  is a pair of operators of the following form
\vskip -13pt
\begin{equation}
\D_g\, \varphi := \D (g\ast \varphi); \; \; \I_f \, \varphi = f\ast \varphi
\end{equation}
satisfying the relations
\begin{eqnarray}
\D_g \, \I_f \, \varphi = \varphi  \; \; \; \forall \varphi \in \mathcal{L}^1_{\mathrm{loc}}(]0,\infty[) \label{x}, \\
\I_f \, \D_g \, \varphi = \varphi  \; \; \; \forall \varphi \in AC([0,\infty[), \label{y}
\end{eqnarray}
where $AC([0,\infty[)$ denotes the set of absolutely continuous functions on $[0,\infty[$.

A sufficient condition for equation~\eqref{x} is the identity
\vskip -10pt
\begin{equation} \label{Sonine}
(g\ast f)(t) = 1 \; \; \; (t \geq 0)
\end{equation}
\vskip -2pt \noindent
because $\D (1\ast \varphi)(t) = \D \int_0^t \varphi(s)\, \dd s = \varphi(t)$ for $t \geq 0 $ and every
$\varphi \in \mathcal{L}^1_{\mathrm{loc}}([0,\infty[)$.

If for a given function $g \in \mathcal{L}^1_{\mathrm{loc}}(]0,\infty[)$ equation~\eqref{Sonine} has a solution $f \in \mathcal{L}^1_{\mathrm{loc}}(]0,\infty[)$, then the solution $f$ is unique in $\mathcal{L}^1_{\mathrm{loc}}([0,\infty[)$, because
$ g\ast F = 0$, $F \in \mathcal{L}^1_{\mathrm{loc}}([0,\infty[)$, implies that $(f \ast g \ast F)(t) = (g\ast f \ast F)(t) =
\int_0^t F(s)\, \dd s = 0$ for all $t > 0$, hence $F(t)\equiv 0$ almost everywhere in $]0,\infty[$.

Equation~\eqref{Sonine} for $f,g \in \mathcal{L}^1_{\mathrm{loc}}(]0,\infty[)$ is known as the Sonine equation \cite{SKM} and its solutions have been studied in \cite{S,K}. The two functions
$g, f$ satisfying the Sonine equation are called a Sonine pair. If $g, f$ is a Sonine pair, then $f, g$ is a Sonine pair and $f$ is
said to be the associated Sonine function of $g$. If $g$ has a Sonine associated function then it is said to be a Sonine function.

We shall show that either of the relations (\ref{x},\ref{y}) is equivalent to equation~\eqref{Sonine}.  
This result has some important consequences for the kernels $g$ and $f$.

\begin{theorem} 
If $f, g \in \mathcal{L}^1_{\mathrm{loc}}$, then
\begin{enumerate}[(i)]
\item $\D_g\, \I_f = \I_f \, \D_g$;
\item $\D_g\, \I_f \, \varphi = \varphi$ for all $\varphi \in \mathcal{C}^1[0,\infty[)$ if and only if $g$ and $f$ satisfy
equation~\eqref{Sonine}.
\end{enumerate}
\end{theorem}
\noindent\textbf{Proof.}\\

The left-hand side of equation~\eqref{x} equals 
$$
\D \int_0^t g(s) \, (f \ast \varphi)(t - s) \, \dd s = g(t) (f \ast \varphi)(0) + \int_0^t g(s)\, \D(f \ast \varphi)(t - s) \, \dd s
$$
The first term equals 0 because $\lim_{t\rightarrow 0} \int_0^t f(\xi)\, \varphi(t - \xi)\, \dd \xi = 0$ by the 
Lebesgue Dominated Convergence Theorem. The second term equals 
\begin{multline} \label{z}
\int_0^t g(s)\, \left[f(t - s)\, \varphi(0) + \left. \int_0^\tau f(\xi) \, \varphi^\prime(\tau - \xi) \, \dd \xi\right|_{\tau=t-s}\right] \, \dd s = \\
\varphi(0) \, (g \ast f)(t) + (g \ast f \ast \varphi^\prime)(t). 
\end{multline} 

We now calculate the left-hand side of equation~\eqref{y}:
\begin{multline*}
(\I_f\, \D_g \, \varphi)(t) = \int_0^t f(t - s)\, \D \int_0^s g(\xi)\, \varphi(s - \xi) \, \dd \xi =
\\
\int_0^t f(t - s) \left[ g(s) \varphi(0) + \int_0^s g(s - \xi)\, \varphi^\prime(\xi)\right] \, \dd \xi = \\
f\ast g \ast \varphi^\prime(t) + \varphi(0)\, (f\ast g)(t) = \D_g \, \I_f \varphi(t)
\end{multline*}
This proves (i).

If equation~\eqref{Sonine} holds then the expression on the right of \eqref{z} equals
$\varphi(0) + (1 \ast \varphi^\prime)(t) \equiv \varphi(t)$, which proves (ii).

\hfill $\Box$
  
\vspace*{-3pt}
\section{Some necessary conditions for a function \break to satisfy the Sonine equation} \label{sec:necessary} 

\setcounter{section}{3} \setcounter{equation}{0}

It is easy to see that a Sonine function is not bounded on any neighborhood of 0. Indeed, if $\vert g \vert \leq M$ on $[0,\varepsilon]$,
for some $M, \varepsilon > 0$, then
$$\left| \int_0^t g(s)\, f(t-s)\, \dd s \right| \leq M\, \int_0^t \vert f(s)\vert\, \dd s \rightarrow 0$$
for $t \rightarrow 0$, which is incompatible with \eqref{Sonine}.

\textit{It follows that the kernel of a GFD cannot be bounded at 0.}

Under an additional assumption it can be proved that the functions $f$, $g$ have the following limits at 0 (\cite{S}):

\begin{theorem}
If $g, f \in \mathcal{L}^1_{\mathrm{loc}}([0,\infty[)$ is a Sonine pair
and \\
($\ast$) the functions $g, f$ are non-increasing over some interval $]0,\varepsilon]$, $\varepsilon > 0$,
then\\
a)
\begin{equation}
\lim_{t\rightarrow 0} g(t) = \infty; \; \; \lim_{t\rightarrow 0} f(t) = \infty,
\end{equation}
b)
\begin{equation}
\lim_{t\rightarrow 0} [t\, g(t)] = 0; \; \;  \lim_{t\rightarrow 0} [t\, f(t)] = 0.
\end{equation}
\end{theorem}

\noindent\textbf{Proof.}\\
Ad a)\\
On account of Assumption ($\ast$) the function $g(t)$ is non-increasing and unbounded on $[0,\varepsilon]$, hence it tends to infinity
for $t \rightarrow 0$.

Similarly for $f(t)$.\\
Ad b)\\
On account of Assumption~($\ast$) in the Sonine equation $$\int_0^t g(s) \, f(t-s) \, \dd s = 1$$
 $g(s) \geq g(t)$ and $f(t - s) \geq f(t)$, hence $g(t)\, f(t)\, t \leq 1$. Since $f(t)$ is unbounded at $t \rightarrow 0$,
$\lim_{t\rightarrow 0} [t\, g(t)] = 0$. Similarly $\lim_{t\rightarrow 0} [t\, f(t)] = 0$.

\hfill$\Box$

{\sc Remark.} Assumption ($\ast$) is satisfied by all the FD definitions proposed so far, including the CF and AB derivatives.\\

If $f, g \in \mathcal{L}^1_{\mathrm{loc}}([0,\infty[)$ and $\tilde{f}, \tilde{g}$ denote their Laplace transforms, then 
equation~\eqref{Sonine} is equivalent to 
\begin{equation} \label{equ}
\tilde{f}(p) \, \tilde{g}(p) = 1/p, \;\; p > 0.
\end{equation}

If the locally integrable functions $g, f$ have regular variation at 0, then
\begin{equation}
g(t) = G(t)\, t^{-\alpha}; \; \; f(t) = F(t)\, t^{-\beta}; \; \; \alpha, \beta < 1,
\end{equation}
where the functions $G$ and $F$ have slow variation at 0 \cite{BGT}. By the Karamata Abelian Theorem
$$ \tilde{g}(p) = H(p) \, p^{\alpha - 1}, \; \;  \tilde{f}(p) = K(p)\, p^{\beta - 1}, $$
where the functions $H, K$ have slow variation at infinity. Equation~\eqref{equ} implies that $\beta = 1 - \alpha$.

If $\alpha \leq 0$ then $\beta \geq 1$, which contradicts local integrability of $f$.

\medskip

We have thus proved that the following theorem.

\begin{theorem}
Suppose that the two locally integrable functions $g$ and $f$ on $[0,\infty[$ have regular variation at 0 with indices $-\alpha$
and $-\beta$ and satisfy the Sonine equation.

If  $0 < \alpha < 1$, then $\beta = 1 - \alpha $ and both $g$ and $f$ have an integrable singularity at 0:
\begin{equation} \label{z}
g(t) = G(t) \, t^{-\alpha}; \; \; f(t) = F(t) \, t^{\alpha -1}, \; \; \; 0 < \alpha < 1,
\end{equation}
where the functions $F, G$ have slow variation at 0.

If the function $f$ is regular at 0 ($\alpha \leq 0$), then its associated Sonine function $g$ is not integrable in any right neighborhood
$[0,\varepsilon]$ of 0, $\varepsilon > 0$.
\end{theorem} 

In particular the slowly varying factors $G, F$ can have a finite value at 0 and be given by convergent power series:
\begin{equation} \label{a}
G(t) = \sum_{n=0}^\infty a_n \, t^n; \; \; F(t) = \sum_{n=0}^\infty b_n \, t^n
\end{equation}
with $a_0 \neq 0$.
Given the coefficients $a_n, n=0, \ldots$ , the coefficients $b_n, n=0, \ldots $ can be determined by a solving sequence of equations
\cite{W,S}. The series for $F$ is assumed convergent, hence in this case we can assert that $g$ satisfies the Sonine equation with the
associated Sonine function $f$ given by equations~\eqref{z} and \eqref{a}.

\section{Completely monotone Sonine pairs} \label{sec:CM}

\setcounter{section}{4} \setcounter{equation}{0}

All the specific examples of generalized fractional derivatives put forward so far as well as the fractional derivatives involve
locally integrable completely monotone (LICM) kernels. Consequently completely monotone solutions of the Sonine equation deserve a separate analysis.

The definition of LICM functions is recalled in Appendix A. 
The LICM function $g$ is said to be a singular LICM (sLICM)
function if $g$ is unbounded near 0 and LICM.

Here is a simple example of a Sonine pair in the sLICM class in closed form:
\begin{multline*}
g(t) = t^{\beta-1} \, E_{\alpha,\beta}\left(-t^\alpha\right); \; \; f(t) = t^{-\beta}/\Gamma(1-\beta) +
t^{\alpha-\beta}/\Gamma(\alpha - \beta + 1),  \\ \beta \geq 0,\; 0 < \alpha\leq \beta,
\end{multline*}
as can be verified using the Laplace transform
$$ \int_0^\infty \exp(-p\, t)\, t^{\beta - 1}\, E_{\alpha,\beta}\left(-t^\alpha\right) \, \dd t = \frac{p^{\alpha-\beta}}{p^\alpha + 1}$$
(\cite{P}, equation~(1.79)). Theorem~2.6 (i-iv) in \cite{GLS} implies that $g$ is sLICM if $\beta \geq 0$, $\alpha \leq \beta$. The
last inequalities entail that $f$ is sLICM.

{\sc Counterexample.}\\
There are Sonine pairs which are not CM, for example
$(\uppi \,t)^{-1/2}\, \cos\left(2\, t^{1/2}\right)$ and $(\uppi\, t)^{-1/2}\, \mathrm{cosh}\left(2\, t^{1/2}\right)$. The corresponding Laplace transforms are $p^{-1/2}\,\\ \exp(-1/p)$ and $p^{-1/2}\, \exp(1/p)$, hence equation~\eqref{equ} is satisfied. The first function
changes sign, hence it is not CM. It follows from Theorem~{\ref{SonineCM}} that the second one is not CM.

The following theorem holds for LICM functions:

\begin{theorem} \label{SonineCM}
Every sLICM function $g: \, ]0,\infty[ \rightarrow \mathbb{R}$ is a Sonine function and the associated Sonine
function $f$ is sLICM.

If $g$ is LICM and regular at 0 then there are a unique positive number $a$ and a unique LICM function $\varphi$ such that
\begin{equation} \label{nS}
a \, g + g \ast \varphi = 1,
\end{equation}
\end{theorem}

\noindent\textbf{Proof.}\\

The functions $f, g$ are a Sonine pair if and only if
their Laplace transforms $\tilde{f}$ and $\tilde{g}$ satisfy equation~\eqref{equ}.

If $g$ is a LICM function, then Bernstein's theorem (Appendix~A) implies that
\vskip -10pt
$$g(t) = \int_{[0,\infty[} \exp(-r\, t) \, \mu(\dd r),$$
where $\mu$ is a Borel measure on $[0,\infty[$ satisfying the inequality~\eqref{ineqmu}.
Hence
$$\tilde{g}(p) = \int_{[0,\infty[} \frac{\mu(\dd r)}{p + r}$$ is a Stieltjes function~\eqref{Stj}.

On account of equation~\eqref{equ} the function
$\tilde{g}$ is not identically 0. Hence $\tilde{g}(p)^{-1}$ is a complete Bernstein function
(Theorem~{\ref{appthm}}). Equation~{\eqref{CBF}} implies that
\vskip -10pt
\begin{equation} \label{c}
\frac{1}{p\, \tilde{g}(p)} = a + \int_{[0,\infty[} \frac{\nu(\dd r)}{p + r},
\end{equation}
where $\nu$ is a Borel measure satisfying the inequality
\begin{equation} \label{ineqnu}
\int_{[0,\infty[} \frac{\nu(\dd r)}{1 + r} < \infty
\end{equation}
and $a \geq 0$.

On account of inequality~\eqref{ineqnu} and the Lebesgue Dominated Convergence Theorem
\vskip -10pt
$$\lim_{p \rightarrow \infty} \frac{1} {p\, \tilde{g}(p)} = a.$$

If $a > 0$ then equation~\eqref{equ} implies that
$$\lim_{t\rightarrow 0} g(t) = \lim_{p\rightarrow \infty} [p\, \tilde{g}(p)] = 1/a < \infty,$$
hence $g(t)$ is not singular at 0. Thus if $g$ is singular at 0 then $a = 0$,
and $f(t) = \int_{[0,\infty[} \exp(-r \,t)\, \nu(\dd r)$ is a LICM solution of the Sonine equation.

If $g$ is regular at 0 then, in view of the monotonicity of $g$,  $0 < g(0) < \infty$, $a = 1/g(0) > 0$ and from \eqref{c}
$$ a\, \tilde{g}(p) + \tilde{g}(p)\, \tilde{\varphi}(p) = 1/p,$$
which is equivalent to equation~{\eqref{nS}}, where $$\varphi(t) := \int_{[0,\infty[} \e^{-r t}\, \nu(\dd r)$$ is a LICM function.

If $f$ is LICM and regular at 0 then it has a positive finite limit $b$ at 0 and
\vskip -10pt
$$\lim_{p \rightarrow \infty} [p \, \tilde{f}(p)] = b.$$
If $f$ is also the solution of the Sonine equation, then it follows from equation~\eqref{equ} that
\vskip -10pt
$$\lim_{p\rightarrow \infty} \tilde{g}(p) = 1/b.$$
Hence
\vskip -10pt
$$\lim_{t \rightarrow 0} \int_0^t g(s)\, \dd s = 1/b > 0$$
which is impossible because $g$ is assumed locally integrable.
Therefore if $f$ is a LICM Sonine function then it is unbounded at 0.

\hfill$\Box$

\vspace{0.15cm}

In view of the uniqueness of the solution $(a, \varphi)$ of equation~\eqref{nS} the inequality $a > 0$ implies that $g$ is not a Sonine function and equations~\eqref{x}, \eqref{y} do not hold for any $f \in \mathcal{L}^1_{\mathrm{loc}}([0,\infty[)$.

\vspace{0.15cm}

{\sc Examples.}\\ 
1. If $g(t) = \exp(-t)$ then $\tilde{g}(p) = 1/(1 + p)$. Note that $1/[p \, \tilde{g}(p)] = 1 + 1/p$ is not the Laplace transform of 
a function, but equation~\eqref{nS} is satisfied with $a = 1$ and $\varphi(t) = 1$.\\
2. If $g(t) = E_\alpha\left(-t^\alpha\right)$, $0 < \alpha \leq 1$, then $\tilde{g}(p) = p^{\alpha-1}/\left(p^\alpha + 1\right)$ and thus
equation~\eqref{nS} holds with $a = 1$ and $\varphi(t) = t^{\alpha -1}/\Gamma(\alpha)$. 

\vspace{0.2cm}

\newpage

\begin{corollary} \label{C}
$(1)$ A LICM function is a Sonine function if and only if it is singular at 0;\\
$(2)$ the associated Sonine function of an sLICM function is an sLICM function.
\end{corollary}

Corollary~{\ref{C}} provides a rich class of generalized fractional derivative and integral kernels. This class of kernels encompasses
fractional derivative kernels as well as two classes of non-fractional GFD with practical applications: distributed-order derivatives \cite{H3} and shifted fractional derivatives \cite{victoria,BlochTorrey}. 

The fractional derivatives of order $\alpha < 1$ are GFD with 
\begin{gather}
g(t) = \theta^{-\alpha}\\
f(t) = \theta^{\alpha-1}
\end{gather}
In this case we shall use the standard notation in Fractional Calculus:
$$\mathrm{D}^\alpha \, \varphi := \mathrm{\D}_g \, \varphi, \; \; \mathrm{I}^\alpha \, \varphi := \mathrm{\I}_f \, \varphi$$
Clearly
\begin{gather}
\widetilde{\D^\alpha \,\varphi}(p) = p^\alpha \, \tilde{\varphi}(p) \label{yz}\\
\widetilde{\I^\alpha \,\varphi}(p) = p^{-\alpha} \, \tilde{\varphi}(p) \label{zz} 
\end{gather}

\section{sLICM kernels with a logarithmic singularity.} \label{log}

Distributed fractional derivatives deserve special attention because they are used in modeling of slow diffusion. Their kernels 
have logarithmic singularities at $t = 0$.

Consider the simplest distributed-order derivative \cite{C,H3,Kochubei2008} with the kernel $w$ 
\begin{equation}
\int_{]0,1]} \D^\alpha \, \varphi(t)\,\dd \alpha = \D (w\ast \varphi)(t),
\end{equation}
where $\mu$ is a bounded non-negative function on $]0,1]$ and
\begin{equation}
w(t) := \int_{]0,1]} \frac{t^{\alpha - 1}}{\Gamma(\alpha)} \, \dd \alpha.
\end{equation}

The kernel $w(t)$ is a superposition of CM functions with non-negative coefficients, hence it is CM.

Furthermore $w$ is integrable over $[0,1]$: $t^{\alpha-1}/\Gamma(\alpha) \leq t^{\alpha-1}$ for $\alpha \in ]0,1]$ and $\int_{]0,1]} t^{\alpha-1} \, \dd \alpha  = 
(t - 1)/[t \ln t]$ for $t > 0$. The majorant is integrable over $[0,1]$. Consequently $w(t)$ is integrable over $[0,1]$ and smooth elsewhere, hence it is locally integrable. 

The family of functions $\alpha \rightarrow t^{1-\alpha}/\Gamma(\alpha)$ is monotone with respect to the parameter $t$ and it tends to infinity
for $t \rightarrow 0$, 
hence $\lim_{t\rightarrow 0} w(t) = \infty$ by a Lebesgue theorem. 

We have thus proved that the function $w$ is sLICM. Consequently it is a Sonine function.

The Laplace transform of $w$ is $\tilde{w}(p) = \int_0^1 p^{-\alpha} \, \dd \alpha = (1 - p)/[p \ln(p)]$.

The associated Sonine function $v$ of $w$ is given by the equation~\eqref{equ}, thus $\tilde{v}(p) = \ln(p)/(1 - p)$.
By a complex contour computation 
$$v(t) = \frac{1}{2\upi \ii} \left[\int_0^\infty \e^{-r t} \frac{\ln\left\{r \,\e^{\ii \upi}\right\}}{1 + r} \, \dd r - 
\int_0^\infty \e^{-r t} \frac{\ln\left\{r \,\e^{-\ii \upi}\right\}}{1 + r} \, \dd r\right] = \int_0^\infty \frac{\e^{-r t}}{1+r} \, 
\dd r.$$

We thus see that $v$ is a superposition of CM functions with non-negative coefficients, hence it is CM. By a similar argument as above $v$ is singular.

We now prove that $v$ is integrable over $[0,1]$, hence locally integrable:
$$\int_0^1 v(t) \, \dd t = \int_0^\infty \frac{\dd s}{1 + s} \int_0^1 \e^{-s t} \, \dd t = \int_0^\infty \frac{\dd s}{s\, (1 + s)} \left(1 - \e^{-s}\right) < \infty.$$ 

Summarizing, $v$ is sLICM in accordance with Theorem~\ref{SonineCM}.

We now show that $w$ has a logarithmic singularity at 0. The asymptotic behavior of $\tilde{w}$ at infinity is $\tilde{w}(p) \sim_\infty
- 1/\ln(p)$. By the Karamata Tauberian Theorem \cite{BGT} the indefinite integral $W(t) := \int_0^t w(t) \, \dd t \sim_0 1/\ln(t)$, hence
$w(t) \sim_0 -1/\left[t \left[ \ln(t)\right]^2\right]$.

The asymptotic behavior of the associated function $v(t)$ will be determined usinf equation~\eqref{equ}:
$\tilde{v}(p) = 1/[p \, \tilde{w}(p)] = \ln(p)/(1 - p) \sim_\infty -p^{-1}\, \ln(p)$. The indefinite integral $V$ of $v$ is non-decreasing and by the Karamata Tauberian theorem has the asymptotic behavior 
$-t \, \ln(t)$ at 0. Hence $v(t) \sim_0 -\ln(t)$.

The asymptotics of more general kernels $w(t) = \int_0^1 t^{\alpha-1} \mu(\alpha)\, \dd \alpha$ is investigated in \cite{Kochubei2008}.

\section{Shifted Fractional Derivatives.} \label{shift}

{\em Shifted Fractional Derivatives} (SFD) \cite{victoria} are defined by the formula 
\begin{equation} \label{xd}
\mathrm{D}^{(\alpha,a)} \, \varphi(t) := \exp(-a\, t) \, \mathrm{D}^\alpha \, \left[ \exp(a \,t)\, \varphi(t) \right],\; \; a > 0
\end{equation}
The Laplace transform of the left-hand side can be readily calculated using equation~\eqref{yz}:
\begin{multline} \label{xy}
\int_0^\infty \exp(-(p + a)\, t) \, \D^\alpha\left(\exp(a\, t) \varphi(t)\right)\, \dd t = \\
\left(\D^\alpha \, \left(\exp(a\, t) \, \varphi(t) \right)\right)\tilde\,(p + a) = p^\alpha \, \tilde{\varphi}(p - a)\vert_{p \rightarrow p + a} 
= \\ (p + a)^\alpha \, \tilde{\varphi}(p)
\end{multline}

Shifted fractional integrals (SFI) are defined by the formula
\begin{equation} \label{xd1}
\mathrm{I}^{(\alpha,a)} \, \varphi(t) := \exp(-a\, t) \, \mathrm{I}^\alpha \, \left[ \exp(a \,t)\, \varphi(t) \right].
\end{equation}
The Laplace transform of the SFI $\mathrm{I}^{(\alpha,a)}$ is given by the equation 
\begin{equation}\label{yy}
\widetilde{\I^{(\alpha,a)} \, \varphi}(p) = (p + a)^{-\alpha}
\end{equation}

From equations~(\ref{xy},\ref{yy}) it follows that 
\begin{equation}
\I^{(\alpha,a)} \, \D^{(\alpha,a)}\, \varphi = \D^{(\alpha,a)} \, \I^{(\alpha,a)} \, \varphi = \varphi
\end{equation}

We shall now show that the shifted fractional derivative is a GFD. Indeed, define the kernel $g_{(\alpha,a)}$ by its Laplace transform 
$\widetilde{g_{(\alpha,a)}}(p) = (p + a)^\alpha/p$. Hence
$$g_{(\alpha,a)}(t) = \frac{1}{2 \ii \uppi} \int_{-\ii \infty}^{\ii \infty} \exp(p \, t)\, (p + a)^\alpha\, \dd p/p$$
Deforming the Bromwich contour to a line running from $-\infty$ below the branching cut $]-\infty,-a]$ and above the cut 
from $-a$ to $-\infty$ we obtain the following formula
\begin{equation} \label{g}
g_{(\alpha,a)}(t) = \exp(-a \, t) \,\theta^{-\alpha}\, \Phi(t,a,\alpha)
\end{equation}
where
\begin{equation}
\Phi(t,a,\alpha) := \frac{1}{\Gamma(\alpha)} \int_0^\infty s^\alpha \, \exp(-s) \, \dd s/(s + a\, t).
\end{equation}
is the superposition of completely monotone functions and therefore it is completely monotone. 
The function $g_{(\alpha,a)}$ is a product of three CM functions, hence it is CM. Since $\Phi(t,a,\alpha) \leq 1$,
the function $g_{(\alpha,a)}$ is singular LICM.

The SFD is thus a GFD with the kernel $g_{(\alpha,a)}$:
\begin{equation}
\mathrm{D}^{(\alpha,a)}\, \varphi = \mathrm{D}_{g_{(\alpha,a)}}\, \varphi
\end{equation}
where $g_{(\alpha,a)}$ is given by equation~\eqref{g}. 

In view of equation~\eqref{equ} the Sonine associated function of $g_{(\alpha,a)}$ is defined by its Laplace transform 
$\tilde{f}_{(\alpha,a)}(p) = (p + a)^{-\alpha}$. Inverting the Laplace transformation we obtain the formula
\begin{equation} \label{f}
f_{(\alpha,a)}(t) = \exp(-a\, t)\, \theta^{\alpha-1}(t)
\end{equation}

This fact has a provides a useful formula for the inverse 
of a SFD:
\begin{equation}
f_{(\alpha,a)} \ast \mathrm{D}^{(\alpha,a)} \, \varphi = \varphi.
\end{equation}

\section{An existence proof for a Sonine pair} \label{sec:exist} 

\setcounter{section}{5} \setcounter{equation}{0}

It may often be difficult to determine the associated Sonine function. We shall consider here a particular generalized fractional
derivative with an sLICM kernel $g$.

We shall prove that
\begin{equation} \label{this}
g(t) := t^{-\alpha}\, \exp(-t^\beta), \; \; 0 < \alpha, \beta < 1
\end{equation}
is a Sonine function and thus it is eligible as a kernel of a generalized fractional derivative.

We begin with noting that $g$  is sLICM. The function $g$  is completely monotone because (1) $t^\beta$, $0 < \beta < 1$, is a Bernstein function; (2) the composition of the completely monotone function $\exp(-x)$ with the Bernstein function $t^\beta$ is CM \cite{SSV}, (3) the product of two completely monotone functions is completely monotone \cite{SSV},
(4) $g(t)$ is integrable on $[0,1]$. Furthermore $g$ is singular at 0.

Theorem~{\ref{SonineCM}} implies that the Sonine equation $g \ast f = 1$ has a unique solution $f$ and $f$ is an
sLICM function.

The function $g$ defined by equation~\eqref{this} is the product of $t^{-\alpha}$ and the function
$G(t) := \exp(-t^\beta)$ slowly varying at 0. Consequently $g$ has regular variation at 0  with the index $-\alpha$. The Karamata
Abelian Theorem \cite{BGT} implies that $\tilde{g}(p) = L(p) \, p^{\alpha-1}$, where $L(p)$ is slowly varying at $\infty$. In view of equation~\eqref{equ}
the Laplace transform of $f$ is $\tilde{f}(p) = L(p)^{-1} \, p^{-\alpha}$. Note that $L(p) \neq 0$ because in view of the inequality $g(t) > 0$
the function $\tilde{g}(p) = \int_0^\infty \e^{-p t}\, g(t)\, \dd t $ is positive for every $p > 0$.
By the Karamata Tauberian Theorem \cite {BGT} $f(t) = F(t) \, t^{\alpha-1}$,
where $F$ is slowly varying at 0.

We now make a working hypothesis that $F(t)$ can be expressed by a convergent power series
\vskip -12pt
$$F(t) =
C\, \left(1 + \sum_{n=1}^\infty a_n\, t^{n \gamma}\right)$$
\vskip -3pt \noindent
for some constants $C, \gamma, a_n, n = 1,2,\ldots$.
The functions $f, g$ are given by the power series
\begin{eqnarray}
g(t) = t^{-\alpha} \, \left[1 - t^\beta + t^{2 \beta}/2! - t^{3 \beta}/3! + \ldots \right], \\
f(t) = C \, t^{\alpha-1}  \, \left[1 + a_1\, t^\gamma + a_2\, t^{2 \gamma}/2! + a_3\, t^{3 \gamma}/3! + \ldots \right], %% . %%
\end{eqnarray}
and their Laplace transforms are given by the formulae
\begin{eqnarray}
\tilde{g}(p) = p^{\alpha-1}\, \Gamma(1-\alpha) \, \left[1 - p^{\alpha-\beta-1}\, \Gamma(1-\beta+\alpha) + \ldots \right],  \\
\tilde{f}(p) = C \, p^{-\alpha} \, \Gamma(\alpha) \, \left[1 + a_1\, p^{-\alpha-\gamma} \, \Gamma(\gamma+\alpha) + \ldots \right].
\end{eqnarray}

We now calculate the left-hand side of the equation $p\, \tilde{f}(p)\, \tilde{g}(p) = 1$ term by term and choose the index $\gamma$ and the coefficients $a_n$ in such a way that all the terms cancel except for the first term, which is equal to 1. We then obtain $\gamma = \beta$, $C = \sin(\uppi \alpha)/\uppi$, $a_1 = \Gamma(\alpha)\, \Gamma(\beta-\alpha+1)/\Gamma(\alpha+\beta) \Gamma(1-\alpha)$. The function $F$ is thus obtained in the form $F(t) = \Phi(t^\beta)$, where the function $\Phi$ is defined by a power series. Convergence of
the power series $\Phi\left(t^\beta\right)$
is ensured because we have proved that the solution $f$ of the Sonine equation for the given Sonine kernel $g$ exists, has the
form $C \, t^{\alpha-1}\, \Phi(t^\beta)$  and is unique.

We have thus proved that $g(t) = t^{-\alpha}\, \exp\left(t^{-\beta}\right)$ is a Sonine kernel and its associated Sonine kernel is a LICM function $f(t) = \Phi\left(t^\beta\right)\, t^{\alpha-1}$, where $\Phi(y)$ is a convergent power series.

\vspace*{-3pt}

\section{Conclusions}  

\setcounter{section}{6} \setcounter{equation}{0}

The analysis of admissible locally integrable kernels of generalized fractional derivatives and integrals has been reduced to the study of the solutions of the Sonine equation.

It has been shown that the kernels of generalized fractional derivatives and integrals cannot be regular at 0. Under some additional technical assumptions they are weakly singular at 0.

The kernels of all the specific generalized fractional derivatives discussed in the literature are completely monotone. In this class of functions
weak singularity of the kernels has been proved without additional technical assumptions. It has also been proved that every singular LICM function is
a Sonine function.

Given the kernel of a generalized fractional derivative the kernel
of the associated generalized fractional integral is in general difficult to calculate numerically. These problems do not appear in
the models with regular kernels because of their simplicity,  which may explain their popularity.

\vspace*{-3pt}

\section*{Acknowledgements}

The author thanks Francesco Mainardi for supplying literature on fractional calculus with regular kernels.

\appendix

\renewcommand{\theequation}{\Alph{section}.\arabic{equation}} %%% added by VK

\section{Some mathematical background used in the paper.}

A function $f$ satisfying the inequalities
$$(-1)^n \, \D^n f(t) \geq 0,  \; \; n = 0,1,\ldots$$
is called a completely monotone (CM) function. A CM function is infinitely differentiable on $]0,\infty[$ and can tend to infinity at 0.

A LICM function is a CM function which is locally integrable (or, equivalently, integrable over $[0,1]$).

A Bernstein function is a non-negative function $f$ such that its derivative is CM. Every Bernstein function is the sum of a primitive of a LICM function and a constant.

\begin{theorem} \label{thm:Bernstein} {\rm (Bernstein)}
A real function $f$ on $]0,\infty[$ is CM if and only if there is a Borel measure $\mu$ on $[0,\infty[$ such that
$$f(t) = \int_{[0,\infty[} \e^{-r t}\, \mu(\dd r), \;\;\;\; t > 0.$$
\end{theorem}

A CM function $f$ is integrable over the interval $[0,1]$ if and only if the measure $\mu$ satisfies inequality~\eqref{ineqmu}
\begin{equation} \label{ineqmu}
\int_{[0,\infty[} \frac{\mu(\dd r)}{1 + r} < \infty,
\end{equation}
(\cite{H1} Appendix~D).
Such CM functions are locally integrable. Locally integrable CM functions are called LICM in this paper.\\

Concerning Stieltjes functions and complete Bernstein functions, see \cite{SSV}. We shall use the following definitions
of such functions: \\

(1) A Stieltjes function $f$ on $[0,\infty[$ is a function which can be expressed by the formula
\vskip -10pt
\begin{equation} \label{Stj}
f(p) = a + \int_{[0,\infty[} (p + r)^{-1} \, \mu(\dd r)
\end{equation}
where $ a \geq 0$ and $\mu$ is a Borel measure satisfying \eqref{ineqmu} (\cite{SSV}, Definition~2.1).

(2) A complete Bernstein function $g$ on $[0,\infty[$ is a function which can be expressed by the formula
\vskip -10pt
\begin{equation} \label{CBF}
f(p) = a\, p  + p \int_{[0,\infty[} (p + r)^{-1} \, \mu(\dd r)
\end{equation}
where $ a \geq 0$ and $\mu$ is a Borel measure satisfying \eqref{ineqmu} (\cite{SSV}, Eq.~(6.5)).\\

In Section~{\ref{sec:CM}} the following theorem (Theorem~7.3  in \cite{SSV}) is used:
\begin{theorem} \label{appthm}
A function $f$ which is not identically 0 is a complete Bernstein function if and only if $1/f$ is a Stieltjes function.
\end{theorem}

%%%%%%%%%%%%%%%%%%%%%%%%%%%%%%%%%%%%%%%%%%%%%%

\bigskip \smallskip

 \it

 \noindent
%(First) Author's full postal address
ul. Bitwy Warszawskiej 1920r 14/52, \\
02-366 Warszawa, POLAND
 \\[4pt]
e-mail: ajhbergen@yahoo.com

\end{document}